\documentclass[12pt]{amsart}

\usepackage[T1]{fontenc}
\usepackage[utf8]{inputenc}
\usepackage[a4paper,margin=2.5cm]{geometry}
\usepackage{amsmath,amssymb}
\usepackage[hidelinks]{hyperref}
\theoremstyle{plain}
\newtheorem{theorem}{Theorem}
\theoremstyle{remark}
\newtheorem{remark}{Remark}

\title[A Chung--Fuchs type theorem]{A Chung--Fuchs type theorem for skew-product dynamical systems}
\author{Xiong Jin}
\address{Department of Mathematics, The University of Manchester, Oxford Road, Manchester M13 9PL, UK}
\email{xiong.jin@manchester.ac.uk}
\date{}

\begin{document}

\begin{abstract}
We prove a Chung--Fuchs type theorem for skew-product dynamical systems with invertible base. Suppose that the Birkhoff averages of a measurable function converge to zero almost surely and that its successive iterates are independent on almost every fibre, with a nonconstant conditional law for the function itself. Then the Birkhoff sums are unbounded above with positive probability and unbounded below with positive probability. If the skew product is ergodic, they oscillate almost surely. As an application, we prove that Mandelbrot cascades acting on the fibre measures of the natural extension of an ergodic symbolic measure are degenerate in the critical case, provided that the weights are not of $\mathcal{Q}_\alpha$-type.
\end{abstract}

\maketitle

\section{Introduction}

Let $(\Omega,\mathcal{A},\mathbb{P})$ be a probability space and let $\{W_n\}_{n\ge 1}$ be a sequence of independent and identically distributed (i.i.d.) real random variables. Write
\[
S=\{S_n=W_1+\cdots +W_n\}_{n\ge 1}
\]
for the associated random walk. We say that $S$ is recurrent if for every $c>0$,
\[
\mathbb{P}(|S_n|\le c \text{ for infinitely many } n)=1.
\]
We recall the classical theorem of Chung and Fuchs \cite{CF51}.
\begin{theorem}[Chung and Fuchs 1951]
If $\mathbb{E}(|W_1|)<\infty$ and $\mathbb{E}(W_1)=0$, then $S$ is recurrent. Furthermore, if $\mathbb{P}(W_1\neq 0)>0$, then $S$ oscillates, that is
\begin{equation}\label{osc}
\mathbb{P}\big(\liminf_{n\to\infty} S_n=-\infty\big)=\mathbb{P}\big(\limsup_{n\to\infty} S_n=+\infty\big)=1.
\end{equation}
\end{theorem}
By the strong law of large numbers, the conditions $\mathbb{E}(|W_1|)<\infty$ and $\mathbb{E}(W_1)=0$ imply
\begin{equation}\label{sll}
\mathbb{P}\Big(\lim_{n\to\infty} \frac{S_n}{n}=0\Big)=1.
\end{equation}
Conversely, \eqref{sll} implies $W_n/n\to0$ almost surely. Independence and the Borel--Cantelli lemma then give $\sum_{n\ge1}\mathbb{P}(|W_1|>n)<\infty$, hence $\mathbb{E}(|W_1|)<\infty$; the strong law now gives $\mathbb{E}(W_1)=0$. Thus, for i.i.d. increments, the recurrence conclusion can equivalently be formulated as
\[
\eqref{sll} \ \Rightarrow \ S \text{ is recurrent}.
\]
For \eqref{osc}, Chung and Fuchs showed that the set of recurrent points is either empty or a closed additive subgroup of $\mathbb{R}$. In the recurrent case, if $\mathbb{P}(W_1\ne0)>0$, this subgroup is nontrivial and is therefore either $\mathbb{R}$ or $a\mathbb{Z}$ for some $a>0$. Either case implies \eqref{osc}.

One may ask for an analogue of the Chung--Fuchs theorem for more general random walks. Let $(X,\mathcal{B},T,\mu)$ be a measure-preserving dynamical system (m.p.d.s.) with $\mu(X)=1$, and let $f:X\to\mathbb{R}$ be measurable. Write
\[
S_f=\{S_nf=f+f\circ T+\cdots +f\circ T^{n-1}\}_{n\ge 1}
\]
for the associated sequence of Birkhoff sums, or generalised random walk. We say that $S_f$ is recurrent if for every $c>0$,
\[
\mu(|S_nf|\le c \text{ for infinitely many } n)=1.
\]
The following extension of the Chung--Fuchs theorem to generalised random walks is essentially due to Dekking \cite{D82}; see also Atkinson \cite{A76} when $f\in L^1(\mu)$.

\begin{theorem}[Dekking 1982]\label{Dekking}
If $\mu$ is ergodic and
\begin{equation}\label{bet}
\mu\Big(\lim_{n\to\infty} \frac{S_nf}{n}=0\Big)=1,
\end{equation}
then $S_f$ is recurrent.
\end{theorem}

For generalised random walks, \eqref{bet} together with $\mu(f\ne0)>0$ need not imply oscillation. A simple counterexample is given by $X=\{x_0,x_1\}$, with $T(x_0)=x_1$, $T(x_1)=x_0$ and $\mu(\{x_0\})=\mu(\{x_1\})=1/2$. Take $f(x_0)=1$ and $f(x_1)=-1$. Then the Birkhoff sums are periodic and bounded, although \eqref{bet} holds.

\medskip

It is natural to ask which independence assumptions, between those of random walks and generalised random walks, ensure oscillation under \eqref{bet}. In \cite{BJ21}, Barral and Jin studied the action of Mandelbrot cascades on ergodic symbolic measures. To establish degeneracy of this action, it was sufficient to show that
\begin{equation}\label{ui}
\mathbb{Q}\Big(\limsup_{n\to\infty}(S_nf+S_nF)=+\infty\Big)=1,
\end{equation}
where $\mathbb{Q}$ is a skew-product extension of an ergodic measure $\mu$, $S_nf$ is a generalised random walk on the base, and $S_nF(x,\cdot)$ is a random walk under the conditional measure $\mathbb{Q}_x$ for $\mu$-almost every $x$. If $(S_nf+S_nF)/n\to c>0$ almost surely, then \eqref{ui} follows immediately. In the cascade setting below, this drift is $h_W-h_\nu$. The critical case, where the drift is zero, requires control of the fluctuations and motivates the following extension of the argument in \cite{BJ21}.

\subsection*{Main result}

Let $(X,\mathcal{B})$ and $(\Omega,\mathcal{A})$ be standard Borel spaces. Let $(X,\mathcal{B},T,\mu)$ be a m.p.d.s. with $\mu(X)=1$, where $T$ is invertible and $T^{-1}$ is measurable. Consider a skew-product m.p.d.s. $(X\times\Omega,\mathcal{F},\hat{T},\mathbb{Q})$, where
\begin{itemize}
\item[-] $\mathcal{F}=\mathcal{B}\otimes\mathcal{A}$;
\item[-] $\hat{T}(x,\omega)=(T(x),g_x(\omega))$, with $(x,\omega)\mapsto g_x(\omega)$ jointly measurable;
\item[-] $\mathbb{Q}\circ\Pi_X^{-1}=\mu$, where $\Pi_X(x,\omega)=x$.
\end{itemize}
Disintegration with respect to $\Pi_X$ gives a measurable family of probability measures $\{\mathbb{Q}_x\}$ on $(\Omega,\mathcal{A})$ such that, for every nonnegative measurable or $\mathbb{Q}$-integrable function $H$,
\[
\int_{X\times\Omega} H(x,\omega)\,\mathbb{Q}(\mathrm{d}x,\mathrm{d}\omega)
=\int_X\int_\Omega H(x,\omega)\,\mathbb{Q}_x(\mathrm{d}\omega)\mu(\mathrm{d}x).
\]
Since $T$ is invertible, the $\hat{T}$-invariance of $\mathbb{Q}$ and uniqueness of disintegration imply
\begin{equation}\label{dis}
\mathbb{Q}_x\circ g_x^{-1}=\mathbb{Q}_{T(x)}
\qquad\text{for }\mu\text{-almost every }x.
\end{equation}
Indeed, the conditional law over $y$ of the push-forward of $\mathbb{Q}$ by $\hat{T}$ is $\mathbb{Q}_{T^{-1}(y)}\circ g_{T^{-1}(y)}^{-1}$. Invariance identifies this law with $\mathbb{Q}_y$. Invertibility of $\hat{T}$ is not required.

For a measurable function $F:X\times\Omega\to\mathbb{R}$ and $n\ge1$, define
\[
\hat{S}_nF=F+F\circ \hat{T}+\cdots+F\circ \hat{T}^{n-1}.
\]
We impose the following assumption:
\begin{itemize}
\item[(a)] for $\mu$-a.e. $x\in X$,
\begin{enumerate}
\item[(i)] $F(x,\cdot)$ is not $\mathbb{Q}_x$-a.s. a constant;
\item[(ii)] the sequence $\{F(x,\cdot),F\circ \hat{T}(x,\cdot),F\circ \hat{T}^2(x,\cdot),\cdots\}$ is an independent sequence of random variables under $\mathbb{Q}_x$.
\end{enumerate}
\end{itemize}
We have the following extension of the Chung--Fuchs theorem.
\begin{theorem}\label{mainthm}
Under the preceding assumptions on the skew product, assume (a). If
\[
\mathbb{Q}\Big(\lim_{n\to\infty} \frac{\hat{S}_nF}{n}=0\Big)=1,
\]
then
\[
\mathbb{Q}\Big(\liminf_{n\to\infty} \hat{S}_nF=-\infty\Big)>0  \text{ and }  \mathbb{Q}\Big(\limsup_{n\to\infty} \hat{S}_nF=+\infty\Big)>0.
\]
\end{theorem}

\begin{remark}
If we further assume that $\mathbb{Q}$ is $\hat{T}$-ergodic, then, since $\{\liminf_{n\to\infty} \hat{S}_nF=-\infty\}$ and $\{\limsup_{n\to\infty} \hat{S}_nF=+\infty\}$ are $\hat{T}$-invariant sets, Theorem \ref{mainthm} implies
\[
\mathbb{Q}\Big( \liminf_{n\to\infty} \hat{S}_nF=-\infty\Big)=\mathbb{Q}\Big( \limsup_{n\to\infty} \hat{S}_nF=+\infty\Big)=1.
\]
\end{remark}

\begin{remark}
Our proof differs from that in \cite{BJ21} in two respects. We use a result of Kesten in place of the stopping time at which $F\circ\hat{T}^n(x,\cdot)$ first becomes positive. We then use an induced transformation to show that a measurable coboundary $u-u\circ T=v$, with $u,v\ge0$ and $u$ finite almost everywhere, satisfies $v=0$ almost everywhere, without assuming that $u$ is integrable.
\end{remark}

\begin{remark}
Assumption (a) holds for a more general class of skew-product measures $\mathbb{Q}$ than the Peyri\`ere measure studied in \cite{BJ21}. In particular, we do not require the law of $F(x,\cdot)$ under $\mathbb{Q}_x$ to be independent of $x$. This allows us to treat the fibre measures of the natural extension of an ergodic symbolic measure, as shown in Section \ref{app}.
\end{remark}

\section{Proof of Theorem \ref{mainthm}}\label{proofs}

We use the filling scheme; see \cite{Der}.

For $a\in\mathbb{R}$, write $a^+=\max\{a,0\}$ and $a^-=\max\{-a,0\}$.

For $n\ge 1$ define
\[
G_n:=\max_{1\le k \le n} \hat{S}_kF.
\]
We have
\[
G_{n+1}^+-G_{n+1}^-=G_{n+1}=F+G_n^+\circ \hat{T}.
\]
Let $G=\lim_{n\to\infty} G_n$. Then
\begin{equation}\label{Gpm}
G^+-G^-=F+G^+\circ \hat{T}.
\end{equation}
By assumption,
\[
\mathbb{Q}\Big(\lim_{n\to\infty} \frac{\hat{S}_nF}{n}=0\Big)=1
\]
so Kesten's result \cite{Ke}, in the form given in \cite[Section 5(c)]{Der}, implies that $\hat{S}_nF$ changes sign infinitely often in the wide sense: it cannot be eventually strictly positive or eventually strictly negative. Although that formulation assumes ergodicity, it applies here to almost every component of the ergodic decomposition of $\mathbb{Q}$, since the almost sure limit is zero on almost every component. Consequently, $G\ge0$ and $G^-=0$ almost surely, and \eqref{Gpm} becomes
\begin{equation}\label{coB}
G^+=F+G^+\circ \hat{T}.
\end{equation}

Assume that $\mathbb{Q}(G^+<\infty)=1$. This implies that for $\mu$-a.e. $x\in X$,
\[
\mathbb{Q}_x(G^+(x,\cdot)<\infty)=1.
\]

Fix $t\in\mathbb{R}\setminus \{0\}$. For $\mu$-a.e. $x\in X$ denote by
\[
\Phi_t(x):=\mathbb{E}_{\mathbb{Q}_{x}}(e^{it G^+(x,\cdot)}) \text{ and } \phi_t(x):=\mathbb{E}_{\mathbb{Q}_{x}}(e^{it F(x,\cdot)}).
\]
For $\mu$-almost every $x$, $G^+\circ\hat{T}(x,\cdot)$ is measurable with respect to the sigma-algebra generated by $\{F\circ\hat{T}^j(x,\cdot):j\ge1\}$. It is therefore independent of $F(x,\cdot)$ under $\mathbb{Q}_x$ by (a)(ii). Together with \eqref{coB} and \eqref{dis}, this gives
\begin{align*}
\Phi_t(x)=&\ \mathbb{E}_{\mathbb{Q}_{x}}(e^{it G^+(x,\cdot)})\\
=&\ \mathbb{E}_{\mathbb{Q}_{x}}(e^{it (F(x,\cdot)+G^+\circ \hat{T}(x,\cdot))})\\
=&\ \phi_t(x)\times \mathbb{E}_{\mathbb{Q}_{x}}(e^{itG^+(T(x),g_x(\cdot))})\\
=&\ \phi_t(x)\times \mathbb{E}_{\mathbb{Q}_{T(x)}}(e^{itG^+(T(x),\cdot)})\\
=&\ \phi_t(x)\times \Phi_t(T(x)),
\end{align*}
Write $u(x)=-\log|\Phi_t(x)|$ and $v(x)=-\log|\phi_t(x)|$, with $-\log0=\infty$. Then $u,v\ge0$ and, $\mu$-almost surely,
\[
u=v+u\circ T.
\]
For each integer $N\ge1$ such that the set
\[
E=\{x\in X:u(x)\le N\}
\]
has positive measure, let $\mu_E=\mu(E)^{-1}\mu|_E$ and let $T_E(x)=T^{\rho_E(x)}(x)$ be the induced transformation, where $\rho_E(x)=\inf\{n\ge1:T^n(x)\in E\}$ is the first return time. By Poincar\'e recurrence, this return time is finite for $\mu_E$-almost every $x$, and $T_E$ preserves $\mu_E$.

Iterating $u=v+u\circ T$ up to the return time gives the coboundary equation
\[
u-u\circ T_E=v_E,
\]
where
\[
v_E(x)=v(x)+v\circ T(x)+\cdots +v\circ T^{\rho_{E}(x)-1}(x).
\]
On $E$, all terms are finite and $0\le v\le v_E\le u\le N$. Integrating and using the $T_E$-invariance of $\mu_E$, we obtain
\[
\int_E v_E\,\mathrm{d}\mu_E
=\int_E u\,\mathrm{d}\mu_E-\int_E u\circ T_E\,\mathrm{d}\mu_E=0.
\]
Nonnegativity implies $v_E=0$ and hence $v=0$ almost everywhere on $E$. Taking the union over $N$ yields
\begin{equation}\label{m0}
\Phi_t(x)\ne0\quad\Longrightarrow\quad |\phi_t(x)|=1
\qquad\text{for }\mu\text{-almost every }x.
\end{equation}

Apply \eqref{m0} simultaneously to the countable sequence $t_n=1/n$. For $\mu$-almost every $x$, the conditional characteristic function $\Phi_t(x)$ is continuous at $0$ and $\Phi_0(x)=1$, since $G^+(x,\cdot)$ is finite $\mathbb{Q}_x$-almost surely. Thus $\Phi_{1/n}(x)\ne0$, and hence $|\phi_{1/n}(x)|=1$, for all sufficiently large $n$.

Fix such an $x$ and let $U,U'$ be independent random variables with the conditional law of $F(x,\cdot)$ under $\mathbb{Q}_x$. For all sufficiently large $n$,
\[
\mathbb{E}\big(e^{i(U-U')/n}\big)=|\phi_{1/n}(x)|^2=1.
\]
Consequently, $e^{i(U-U')/n}=1$ almost surely for all these $n$, on a common event of probability one. Since $(U-U')/n\to0$, this forces $U=U'$ almost surely. Independence then implies that their common law is a point mass, contradicting (a)(i).

We conclude that $\mathbb{Q}(G^+<\infty)<1$. Since each finite partial sum is real, an unbounded supremum is equivalent to an infinite limsup, so
\[
\mathbb{Q}\Big(\limsup_{n\to\infty} \hat{S}_nF=+\infty\Big)>0.
\]
Applying the same argument to $-F$, we also obtain
\[
\mathbb{Q}\Big(\liminf_{n\to\infty} \hat{S}_nF=-\infty\Big)>0.
\]

\section{Application}\label{app}

We apply Theorem \ref{mainthm} to show that Mandelbrot cascades acting on fibre measures of the natural extension of an ergodic symbolic measure are degenerate in the critical case, when the weights are not of $\mathcal{Q}_\alpha$-type.

\medskip

Let $b\ge2$ and let $X=X^-\times X^+$ be the two-sided symbolic space over $\{1,\ldots,b\}$, where
\[
X^+=\{1,\ldots,b\}^{\mathbb{N}},\qquad
X^-=\{1,\ldots,b\}^{\{\ldots,-1,0\}}.
\]
Let $\Sigma=\bigcup_{n\ge1}\{1,\ldots,b\}^n$ be the set of nonempty finite words, and write $|u|$ for the length of $u\in\Sigma$. For words $u,v$ of lengths $m,n$, respectively, define
\[
[u,v]=\{x\in X:x_{-m+1}\cdots x_0=u,\ x_1\cdots x_n=v\}.
\]
Write $[v]$ for the cylinder in $X^+$ specified by $v$. Let $\mathcal{B}$ and $\mathcal{B}(X^\pm)$ denote the Borel sigma-algebras of the corresponding product topologies. On $X^+$ we use the standard ultrametric
\[
d(x^+,y^+)=b^{-|x^+\wedge y^+|}\quad(x^+\ne y^+),
\]
where $|x^+\wedge y^+|$ is the length of the longest common prefix.

Let $\nu$ be a shift-invariant ergodic probability measure on $X^+$ and let $\mu$ be its natural extension to $X$. We use $\sigma$ for the left shift on either space. On $X$, the shift is invertible and $\mu$ is ergodic.

For $x\in X$, write $x=x^-x^+$ for its past and future. Let $\pi(x)=x^-$ and $\mu^-=\mu\circ\pi^{-1}$. For $\mu^-$-almost every $x^-$, let $\nu_{x^-}$ be the conditional distribution of the future given the past. For every nonnegative measurable or $\mu$-integrable function $f$,
\begin{equation}\label{eq0}
\int_X f(x)\,\mu(\mathrm{d}x)
=\int_{X^-}\int_{X^+}f(x^-x^+)\,\nu_{x^-}(\mathrm{d}x^+)\mu^-(\mathrm{d}x^-).
\end{equation}

Let $(\Omega_0,\mathcal{A}_0,\mathbb{P}_0)$ be a standard Borel probability space and let $W$ be a nonnegative random variable on this space with $\mathbb{E}_{\mathbb{P}_0}(W)=1$. Define
\[
(\Omega,\mathcal{A},\mathbb{P})
=\bigotimes_{u\in\Sigma}(\Omega_0,\mathcal{A}_0,\mathbb{P}_0),
\qquad W_u(\omega)=W(\omega_u).
\]
Then $\{W_u\}_{u\in\Sigma}$ is an i.i.d. family. For $n\ge1$, put
\[
V_n(x,\omega)=W_{x^+|_n}(\omega),
\]
where $x^+|_n=x_1\cdots x_n$.

We construct the Peyri\`ere measure on the full product sigma-algebra $\mathcal{F}=\mathcal{B}\otimes\mathcal{A}$. Let
\[
\mathbb{P}_0^W(\mathrm{d}\omega_0)=W(\omega_0)\mathbb{P}_0(\mathrm{d}\omega_0).
\]
For each $x\in X$, define the product probability measure
\[
\mathbb{Q}_x=\bigotimes_{u\in\Sigma}\rho_{u,x},\qquad
\rho_{u,x}=
\begin{cases}
\mathbb{P}_0^W,&u=x^+|_{|u|},\\
\mathbb{P}_0,&\text{otherwise}.
\end{cases}
\]
Thus the coordinates on the sampled branch have the size-biased law, while all other coordinates retain their original law. This is a measurable probability kernel: its values on cylinder sets depend measurably on finitely many symbols of $x$. Set
\[
\mathbb{Q}(\mathrm{d}x,\mathrm{d}\omega)
=\mu(\mathrm{d}x)\mathbb{Q}_x(\mathrm{d}\omega).
\]
For $n\ge1$, let $\mathcal{A}_n=\sigma\{\omega\mapsto\omega_u:1\le|u|\le n\}$ and $\mathcal{F}_n=\mathcal{B}\otimes\mathcal{A}_n$. Then $\mathcal{F}=\sigma(\bigcup_{n\ge1}\mathcal{F}_n)$, and for every nonnegative $\mathcal{F}_n$-measurable or $\mathbb{Q}$-integrable $\mathcal{F}_n$-measurable function $f$,
\begin{equation}\label{eq1}
\int_{X\times\Omega}f\,\mathrm{d}\mathbb{Q}
=\int_{X\times\Omega}f(x,\omega)\prod_{j=1}^nV_j(x,\omega)\,\mu(\mathrm{d}x)\mathbb{P}(\mathrm{d}\omega).
\end{equation}
In particular, this construction uses only the finite-level densities and does not require non-degeneracy of the limiting cascades.

\medskip

Define the measurable skew-product transformation
\[
\hat{T}(x,\omega)=(\sigma(x),g_x(\omega)),\qquad
g_x(\omega)=(\omega_{x_1u})_{u\in\Sigma}.
\]
The product construction gives $\mathbb{Q}_x\circ g_x^{-1}=\mathbb{Q}_{\sigma(x)}$, because a coordinate $x_1u$ lies on the branch $x^+$ exactly when $u$ lies on the branch of $\sigma(x)$. Together with the $\sigma$-invariance of $\mu$, this proves the $\hat{T}$-invariance of $\mathbb{Q}$.

We also verify ergodicity directly. Let $f,g$ be bounded $\mathcal{F}_m$-measurable functions, and write
\[
\bar f(x)=\int_\Omega f(x,\omega)\,\mathbb{Q}_x(\mathrm{d}\omega),
\qquad
\bar g(x)=\int_\Omega g(x,\omega)\,\mathbb{Q}_x(\mathrm{d}\omega).
\]
For $k>m$, the environment coordinates used by $g$ have depth at most $m$, whereas those used by $f\circ\hat{T}^k$ have depths between $k+1$ and $k+m$. Conditional on $x$, these families are independent under $\mathbb{Q}_x$. Hence
\[
\int_{X\times\Omega}(f\circ\hat{T}^k)g\,\mathrm{d}\mathbb{Q}
=\int_X(\bar f\circ\sigma^k)\bar g\,\mathrm{d}\mu.
\]
Ergodicity of $\mu$ therefore yields
\[
\lim_{N\to\infty}\frac1N\sum_{k=0}^{N-1}
\int_{X\times\Omega}(f\circ\hat{T}^k)g\,\mathrm{d}\mathbb{Q}
=\left(\int f\,\mathrm{d}\mathbb{Q}\right)
 \left(\int g\,\mathrm{d}\mathbb{Q}\right).
\]
Since the sigma-algebras $\mathcal{F}_m$ increase to $\mathcal{F}$, approximation in $L^2(\mathbb{Q})$ extends this identity to all bounded measurable $f,g$. Applying it to the indicator of an invariant set proves ergodicity of $\mathbb{Q}$. This argument applies in particular at criticality.

\medskip

Write $h_\nu$ for the measure theoretic entropy of $\nu$ and
\[
h_W=\mathbb{E}_{\mathbb{P}_0}(W\log W),
\]
with $0\log0=0$. Assume from now on that
\[
h_W=h_\nu<\infty,
\]
which is the critical case. Define
\[
F(x,\omega)=\log\nu_{x^-}([x_1])+\log W_{x^+|_1}(\omega)
\]
where both logarithms are finite, and set $F=0$ on the $\mathbb{Q}$-null exceptional set. The conditional entropy formula for the natural extension gives
\[
-\int_X\log\nu_{x^-}([x_1])\,\mu(\mathrm{d}x)=h_\nu\le\log b;
\]
see also the proof of \cite[Lemma 3.1]{BrJ22}. Moreover,
\[
\int_{X\times\Omega}|\log W_{x^+|_1}|\,\mathrm{d}\mathbb{Q}
=\mathbb{E}_{\mathbb{P}_0}(W|\log W|)<\infty,
\]
since $h_W<\infty$ and $-w\log w\le1/e$ for $0<w\le1$. Thus $F\in L^1(\mathbb{Q})$, and
\begin{align*}
\int_{X\times\Omega}F\,\mathrm{d}\mathbb{Q}
&=\int_X\log\nu_{x^-}([x_1])\,\mu(\mathrm{d}x)
  +\mathbb{E}_{\mathbb{P}_0}(W\log W)\\
&=-h_\nu+h_W=0.
\end{align*}
By Birkhoff's ergodic theorem and ergodicity of $\mathbb{Q}$,
\[
\mathbb{Q}\Big(\lim_{n\to\infty}\frac{\hat{S}_nF}{n}=0\Big)=1.
\]

For every $n\ge1$, almost surely,
\[
F\circ\hat{T}^{n-1}(x,\omega)
=\log\nu_{x^-x_1\cdots x_{n-1}}([x_n])
 +\log W_{x^+|_n}(\omega),
\]
where $x^-x_1\cdots x_{n-1}$ denotes the past of $\sigma^{n-1}(x)$. For $\mu$-almost every fixed $x$, the first term is deterministic and the path weights in the second term are independent under $\mathbb{Q}_x$, with common size-biased law
\[
w\,\mathbb{P}_0(W\in\mathrm{d}w).
\]
This verifies (a)(ii).

Assume also that there is no $\alpha\ge0$ for which
\[
W=\begin{cases}
e^\alpha,&\text{with probability }e^{-\alpha},\\
0,&\text{with probability }1-e^{-\alpha}.
\end{cases}
\]
In the terminology of \cite{BJ21}, the cascade action $\mathcal{Q}$ is not of $\mathcal{Q}_\alpha$-type. Equivalently, the size-biased law of $W$ is not a point mass. Consequently, $F(x,\cdot)$ is nonconstant under $\mathbb{Q}_x$ for $\mu$-almost every $x$, verifying (a)(i). Theorem \ref{mainthm}, together with ergodicity, now gives
\begin{equation}\label{limsup}
\mathbb{Q}\Big(\limsup_{n\to\infty}\hat{S}_nF=+\infty\Big)=1.
\end{equation}

The chain rule for conditional probabilities gives, for $\mu$-almost every $x$ and every $n\ge1$,
\[
\nu_{x^-}([x_1\cdots x_n])
=\prod_{j=1}^n\nu_{x^-x_1\cdots x_{j-1}}([x_j]);
\]
see the proof of \cite[Lemma 3.1]{BrJ22}. The term for $j=1$ is $\nu_{x^-}([x_1])$. Thus, for $\mathbb{Q}$-almost every $(x,\omega)$ and every $n\ge1$,
\begin{equation}\label{snf}
\hat{S}_nF(x,\omega)
=\log\left(\nu_{x^-}([x_1\cdots x_n])
\prod_{j=1}^nW_{x^+|_j}(\omega)\right).
\end{equation}
Combining \eqref{eq0} and \eqref{eq1}, for every nonnegative $\mathcal{F}_n$-measurable or $\mathbb{Q}$-integrable $\mathcal{F}_n$-measurable function $f$, we have
\begin{equation}\label{eq2}
\begin{split}
\int_{X\times\Omega}f\,\mathrm{d}\mathbb{Q}
={}&\int_{X^-}\int_\Omega\int_{X^+}
 f(x^-x^+,\omega)\prod_{j=1}^nW_{x^+|_j}(\omega)\\
&\hspace{18mm}\nu_{x^-}(\mathrm{d}x^+)\mathbb{P}(\mathrm{d}\omega)\mu^-(\mathrm{d}x^-).
\end{split}
\end{equation}
In particular, the conditional probability
\[
\mathbb{Q}^{x^-}(\mathrm{d}x^+,\mathrm{d}\omega)
=\nu_{x^-}(\mathrm{d}x^+)\mathbb{Q}_{x^-x^+}(\mathrm{d}\omega)
\]
is the Peyri\`ere measure associated with the cascade acting on $\nu_{x^-}$. Disintegrating \eqref{limsup} over $\mu^-$ and using \eqref{snf}, we obtain, for $\mu^-$-almost every $x^-$,
\[
\mathbb{Q}^{x^-}\left(
\sup_{n\ge1}\nu_{x^-}([x_1\cdots x_n])
\prod_{j=1}^nW_{x^+|_j}=\infty\right)=1.
\]
By \cite[Theorem 2.1(i)]{BJ21}, $\mathcal{Q}\cdot\nu_{x^-}=0$ $\mathbb{P}$-almost surely for $\mu^-$-almost every $x^-$. More precisely, for these pasts, the almost sure weak limit
\[
\widetilde{\nu}_{x^-}^{(\omega)}(\mathrm{d}x^+)
=\lim_{n\to\infty}\left(\prod_{j=1}^nW_{x^+|_j}(\omega)\right)
\nu_{x^-}(\mathrm{d}x^+)
\]
is the zero measure.

Finally, by \cite[Lemma 3.1]{BrJ22}, $\nu$ and $\mu^-$-almost every $\nu_{x^-}$ are exact-dimensional for the metric $d$, with dimension $h_\nu/\log b$. Hence the disintegration
\[
\nu=\int_{X^-}\nu_{x^-}\,\mu^-(\mathrm{d}x^-)
\]
expresses $\nu$ as a mixture of $\mathcal{Q}$-singular probability measures having the same exact dimension as $\nu$.
One can nevertheless construct an ergodic measure $\nu$ and weights not of $\mathcal{Q}_\alpha$-type such that $h_\nu=h_W$ and $\mathbb{E}_{\mathbb{P}}(\mathcal{Q}\cdot\nu)=\nu$, although $\mathcal{Q}\cdot\nu_{x^-}=0$ $\mathbb{P}$-almost surely for $\mu^-$-almost every $x^-$. Thus fibrewise degeneracy need not imply degeneracy of the cascade acting on the original measure.

\end{document}